# *New Horizons On Fuzzy Magic Graph*

*Jamil, R. N.*[1,a]*, Rehman, M. A.*[1,b]*, Javaid, M.*[2,c]

*University of Management and Technology , c-II Johar Town Lahore Pakistan* [1]

*Department of Mathematical Sciences University of Science and Technology of China, Hefei, China*[2]

*email: noshadjamil@yahoo.com*[a]*, aziz3037@yahoo.com*[b]*,  javaidmath@gmail.com*[c]

**Abstract**: *In this paper, we define the concept of fuzzy magic graphs. A fuzzy graph $G = \{\alpha, \beta\}$ is said to be a fuzzy magic graph if there exist two injective functions $\alpha: V \to [0; 1]$ and $\beta : V \times V \to [0; 1]$ such that $\beta(uv) < \alpha(u) + \alpha(v)$ and $\alpha(u) + \beta(uv) + \alpha(v) = m(G)$ for all $u, v \in V(G)$, where $m(G) \in [0; 1]$ is a fuzzy magic constant. Moreover, we investigate some families of fuzzy graphs like fuzzy paths, fuzzy stars and fuzzy cycles which are fuzzy magic graphs.*

1. Introduction

Fuzzy sets introduced by Zadeh in 1965 [1]. Fuzzy set theory is a very useful concept for uncertain situation in real life problems. While graphs theory first time used by Leonard Euler in 1736 to solve well known "The K¨onigsberg Bridge Problem". Since then, graph theory used in different field of science and technology such that in chemistry to investigation the development of chemical bonds [2], decision making in the manufacturing engineering [3], it also help to better understanding the biological sciences [4], graph theory assist the scientist to examine the some problems in statistical physics [5], in 21th century Computer Science is used in almost every field but without graph theory it is not possible [6], graph theory also used in linguistics, social science etc. [7,8].  Some researcher independently introduced the concepts of fuzzy graph theory by merging the Fuzzy set and graph theory [9, 10]. There are many problems which may be solved by fuzzy graph theory. Here we discuss an example related industrial management, let us consider a company has five departments say $\mathbf{D; D_1; D_2; D_3}$ **and** $\boldsymbol{D_4}$**,** **"D"** coordinates with all other departments. Company gets a project and wants to complete it within a specific time period. Due to this company wishes to allocate work load according to individuals (in a department) working abilities and size of departments. Due to difficulty arise in estimation of human working ability based on qualities, so as well for whole department. But in fuzzy conditions, any person or department working ability always lies within$[0\ \ 1]$. Hence fuzzy magic graph model helps the company to allocate the work load to complete the task within given time frame. In this model boxes represent the departmental work load, edges between boxes shows the share work load between coordinator and other departments.

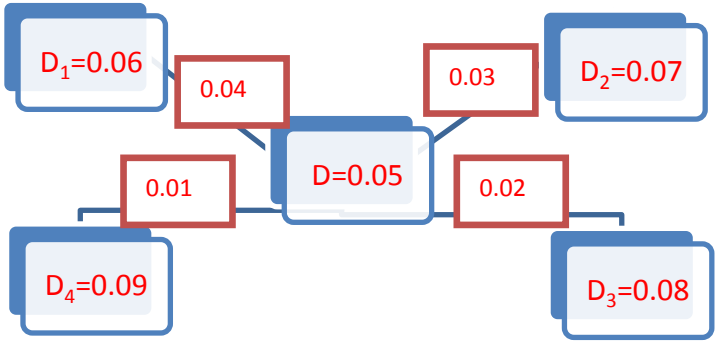

Following table represents % work done by each department and sharing work with D to complete task within specific time period.

| Dep. | D | $D_1$ | $D_2$ | $D_3$ | $D_4$ |
|---|---|---|---|---|---|
| D.'s Work | 33.33% | 40% | 46.67% | 53.33% | 60% |
| Sharing | | $E_1$ | $E_2$ | $E_3$ | $E_4$ |
| S. 's work | | 26.67% | 20% | 13.33% | 6.67% |

Total work done by each department with sharing and coordinator must be 100% for completion of task = $D+D_i + E_i$ ;   $i = 1,2,3,4$

In this paper section 1, include discussion about application of graph theory in different subjects. Section 2 contains some basic definitions while section 3 consists of some new definitions and theorems while the last section gives the conclusion of the whole discussion. All graphs used in this article are connected and finite.

## 2 Preliminaries

Firstly we talk about some definition and results detail can be seen in [9 - 15].

***Definition 2.1:*** Let us consider $X$ be a universal set, then the set $A$ is called a fuzzy set in $X$, if $A = \{(x, f_A(x))\}$, where $x \in X$ and $f_A(x): X \to [0, 1]$ is a function which is known as a membership function e. g. if $X = \{0,1,2,3,4\}$, then $A_1 = \{(0,0.02), (1,0.03), (2,0.05, (3,0.06), (4,0.08)\}$, $A_2 = \{(0,0.01), (1,0.03), (2,0.05, (3,0.07), (4,0.09)\}$ and $A_3 = \{(0,0.00), (1,0.02), (2,0.04, (3,0.06), (4,0.08)\}$ are three fuzzy sets in $X$.

***Definition 2.2:*** If $X$ and $Y$ are two sets and $\mu$ be a fuzzy set of $X \times Y$ then $\mu$ is called fuzzy relation from $X$ into $Y$.

***Definition 2.3:*** A fuzzy relation is said to be symmetric if $\mu(x,y) = \mu(y,x) \; \forall \; x, y \in A$

***Definition 2.4:*** A graph $G(V, E)$ is a pair of vertex-set $V(G)$ and edge-set $E(G)$ such that $E(G) \epsilon \; V(G) \times V(G)$. In addition, $\mathbf{v} = |V(G)|$ and $\mathbf{e} = |E(G)|$ are order and size of the graph $\mathbf{G(V, E)}$.

***Definition 2.5:*** A fuzzy graph $\mathbf{G} = \{\tau, \varphi\}$ with vertex-set $\mathbf{V(G)}$ is a pair of two injective functions $\tau: \mathbf{V(G)} \to [\mathbf{0}, \mathbf{1}]$ and $\varphi: \mathbf{V(G)} \times \mathbf{V(G)} \to [\mathbf{0}, \mathbf{1}]$ such that $\varphi(uv) < \tau(u) \times \tau(V), \forall \; u, v \epsilon \; V(G)$, where $\tau(\mathbf{u})$ and $\tau(\mathbf{v})$ are called membership values of the vertices $\mathbf{u}$ and $\mathbf{v}$ respectively, and µ(uv) iscalled membership value of the edge uv.

***Definition 2.6:*** If $\mathbf{G} = \{\tau, \varphi\}$ be a fuzzy graph then the degree of a vertex $\mathbf{u}$ is $\mathbf{d_G(u)} = \sum_{\mathbf{u} \neq \mathbf{v}} \varphi \; (\mathbf{uv} = \mathbf{e})$, where $\varphi(\mathbf{e}) > 0, \forall \; e \epsilon E$ and $\varphi(\mathbf{e}) = \mathbf{0}, \forall \mathbf{e} \notin E$. $\boldsymbol{\delta}(\mathbf{G}) = \wedge \{\mathbf{d_G(u)}/\mathbf{u} \in \mathbf{V}\}$ is called minimum degree of $\mathbf{G}$ while maximum degree of $\mathbf{G}$ is defined as $\boldsymbol{\Delta}(\mathbf{G}) = \vee \{\mathbf{d_G(u)}/\mathbf{u} \in \mathbf{V}\}$.

Rosenfeld, Bhutani, Battou and Bhattacharya [9,10,13] defied the concept of fuzzy graph for different families of graphs like:

**Definition 2.7:** A path $P_n$ with set of vertices $V(P_n) = \{v_i : 1 \leq i \leq n+1\}$ and set of edges $E(P_n) = \{v_i v_{i+1} : 1 \leq i \leq n+1\}$ is said to be fuzzy path if $\varphi(v_i v_{i+1}) > 0, \forall\, 1 \leq i \leq n$, where $n$ is called the length of the path. A path $P_n$ is called a cycle $C_n$ with length $n$ if $v_1 = v_{n+1}$ and $n \geq 3$.

**Definition 2.8:** A star $S_{1,n} = \{vu_i / 1 \leq i \leq n\}$ is said to be fuzzy star if $\varphi(vu_i) > 0, \forall\, 1 \leq i \leq n$ where $u_i \in U$ and $v \epsilon V$ are sets of vertices such that $|V| = 1$ and $|U| = n \geq 2$.

Nagoor et al. [16] defined the concept of fuzzy labeling as follow:

**Definition 2.9:** A fuzzy graph $G = \{\tau, \varphi\}$ is said to be fuzzy labeling graph if membership values of all vertices and edges are different.

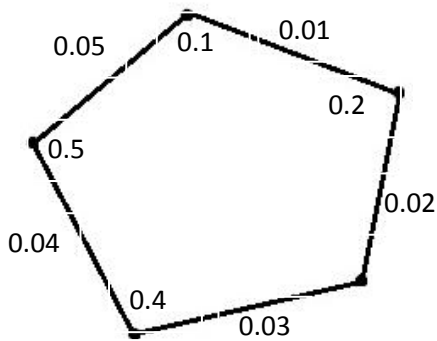

Figure 2. A Fuzzy Labeling Graph

## Main work for discussion

**Definition:** - A graph G= $(\alpha, \beta)$ is said to be a fuzzy magic graph if there exist two injective functions $\alpha : V \to [0\ 1]$ and $\beta : V \times V \to [0\ 1]$ with restricted the following conditions

1. $\beta(uv) < \alpha(u) + \alpha(v)$
2. $\alpha(u) + \beta(uv) + \alpha(v) = m(G) \leq 1$,

Where $m(G)$ is a real constant $\forall\, u, v \in G$

**Definition:** A path "$P_n$" is said to be fuzzy path if $\beta(e_i) > 0$ for all $1 \leq i \leq n$, where $n$ is the length of "$P_n$".

**Definition:** "$S_{1,n}$" is said to be fuzzy star if $\beta(vu_i) > 0$ for all $1 \leq i \leq n$. Where $v \in V$ and $u_i \in U$ are set of vertices such that $O(V) = 1$ and $O(U) = n \geq 2$.

**Theorem 1: A fuzzy Path "$P_n$" is a fuzzy magic graph.**

**Proof:** Let $P_n$ is a fuzzy path with length $n \geq 1$ defined on set of vertices $V$.

**Case 1:** If "$P_n$" has odd length then fuzzy magic labeling as follow:

1. Labeling for odd vertices $\alpha(v_{2i-1}) = (2n + 2 - i)d$      where $1 \leq i \leq \frac{n+1}{2}$
2. Labeling for even vertices $\alpha(v_{2i}) = \alpha(v_n) - id$      where $1 \leq i \leq \frac{n+1}{2}$
3. Labeling for edges $\beta(v_{n+1-i}v_{n+2-i}) = \alpha(v_{n+1}) - id$      where $1 \leq i \leq n$

Where "$d$" is define as follow

$$d = \begin{cases} 10^{-1} & n = 1 \\ 10^{-2} & 3 \leq n \leq 27 \\ 10^{-3} & 29 \leq n \leq 285 \\ 10^{-(i+4)} & 285 \times 10^i < n < 285 \times 10^{i+1} \ \ i = 0,1,2 \ldots \end{cases}$$

**Case 2:** If "$P_n$" has even length then fuzzy magic labeling as follow:

To find the magic constant $m(G)$ one can use

4. $m(G) = \alpha(u) + \beta(uv) + \alpha(v)$      where $1 \leq i \leq \frac{n}{2} + 1$
5. $\phantom{m(G)} = (2n + 2 - i)d + \alpha(v_{n+1}) - id + \alpha(v_n) - id$      where $1 \leq i \leq \frac{n}{2}$
6. $\phantom{m(G)} = (2n + 2 - 3i)d + \alpha(v_{n+1}) + \alpha(v_n)$      where $1 \leq i \leq n$

Where "$d$" is define as follow

$$d = \begin{cases} 10^{-1} & n = 2 \\ 10^{-2} & 4 \leq n \leq 28 \\ 10^{-3} & 30 \leq n \leq 284 \\ 10^{-(i+4)} & 284 \times 10^i < n < 284 \times 10^{i+1} \ \ i = 0,1,2 \ldots \end{cases}$$

$m(G)$ can find as Case 1

## Theorem 2: A fuzzy star "$S_{1,n}$" is a fuzzy magic graph.

**Proof:** Let "$S_{1,n}$" is a fuzzy star with $\beta(vu_i) > 0$ for all $1 \leq i \leq n$. Where $v \in V$ and $u_i \in U$ are set of vertices such that $O(V) = 1$ and $(U) = n \geq 2$.

**Fuzzy magic labeling for "$S_{1,n}$" is defined as follow:**

1. Labeling for vertex "$v$" is $\alpha(v) = (n + 1)d$      where $v \in V$
2. Labeling for vertices $\alpha(u_i) = \alpha(v) + id$      where $1 \leq i \leq n$
3. Labeling for edges $\beta(vu_i) = \alpha(v) - id$      where $1 \leq i \leq n$

If $a = 3\sum_{j=0}^{l-1}(10^j - 1)$, $b = 3\sum_{j=0}^{l-1}(10^j)$ and $l =$ Total number of digits in $n$ (i.e. if $1 \leq n \leq 9 \Rightarrow l = 1$ or if $10 \leq n \leq 99 \Rightarrow l = 2$ and so on), then

$$d = \begin{cases} 10^{-1} & n = 2 \\ 10^{-2} & 3 < n \leq 9 \\ 10^{-l} & 10^{l-1} \leq n \leq a \\ 10^{-(l+1)} & b \leq n \leq 10^l - 1 \end{cases}$$

Thus the fuzzy magic constant is

$$\pi(G) = \alpha(u) + \beta(e) + \alpha(v)$$
$$= (n+1)d + \alpha(v) - id + \alpha(v) + id$$
$$= (n+1)d + (n+1)d + (n+1)d = 3(n+1)d$$

**Theorem 3:** A fuzzy Cycle "$C_n$" is a fuzzy magic graph. where $n \equiv 1(mode 2)$

**Proof:** Let "$C_n$" is a fuzzy Cycle. Consider $V = \{v_1, v_2, v_3, \ldots, v_n\}$ and $E = \{v_1v_2, v_2v_3, v_3v_4, \ldots, v_{n-1}v_n, v_nv_1\}$ be sets of vertices and edges, where $|V| = 2m+1$, $\forall m \in N$

Fuzzy magic labeling for "$C_n$" is defined as follow:

1. $\beta(e) = \begin{cases} \beta(v_i v_{i+1}) = id \\ \beta(v_1 v_n) = nd \end{cases}$  where $\begin{cases} 1 \leq i \leq n-1 \\ i = n \end{cases}$

2. Labeling for odd vertices $\alpha(v_{(n+2)-2i}) = \beta(v_n v_1) + id$  where $1 \leq i \leq \frac{n+1}{2}$

3. Labeling for even vertices $\alpha(v_{(n+1)-2i}) = \beta(v_1) + id$  where $1 \leq i \leq \frac{n-1}{2}$

Where "$d$" is define as follow

$$d = \begin{cases} 10^{-2} & 3 \leq n \leq 27 \\ 10^{-3} & 29 \leq n \leq 287 \\ 10^{-3} & 289 \leq n \leq 2849 \\ 10^{-(i+4)} & 285 \times 10^i + 1 < n < 285 \times 10^{i+1} - 1 \quad i = 1,2 \ldots \end{cases}$$

Thus the fuzzy magic constant is
$$\pi(G) = \alpha(v_i) + \beta(v_i v_{i+1}) + \alpha(v_{i+1}), \text{ for } 1 \leq i \leq n-1$$